\documentclass[11pt]{article} 

\usepackage{setspace}
\usepackage{color}
 \usepackage{amsfonts,amssymb,latexsym,bm,amsmath,amsfonts}
 \usepackage[latin1]{inputenc}  

\usepackage{listings}
\usepackage{xcolor}

\definecolor{codegreen}{rgb}{0,0.6,0}
\definecolor{codegray}{rgb}{0.5,0.5,0.5}
\definecolor{codepurple}{rgb}{0.58,0,0.82}
\definecolor{backcolour}{rgb}{0.95,0.95,0.92}

\lstdefinestyle{mystyle}{
    backgroundcolor=\color{white},   
    commentstyle=\color{codegreen},
    keywordstyle=\color{black},
    numberstyle=\tiny\color{codegray},
    stringstyle=\color{codepurple},
    basicstyle=\ttfamily\footnotesize,
    breakatwhitespace=false,         
    breaklines=true,    captionpos=b,                    
    keepspaces=true, numbers=left,                    
    numbersep=5pt, showspaces=false,                
    showstringspaces=false, showtabs=false, tabsize=2}
\newtheorem{innercustomthm}{Theorem}
\newenvironment{customthm}[1]
  {\renewcommand\theinnercustomthm{#1}\innercustomthm}
  {\endinnercustomthm}

\makeindex    
\newtheorem{thm}{Theorem}

\newtheorem{lemma}{Lemma}
\newtheorem{corollary}{Corollary}
\newtheorem{proposition}{Proposition}
\newtheorem{exam}{Example}

\def\N{\mathbb{N}}
\def\L{{\cal L}}

\def\p{\hbox{\rm\setbox1=\hbox{I}\copy1\kern-.45\wd1 P}}

\def\pr{\noindent{\bf Proof}\ }

\def\l{\lambda}

\def\d{d_{_{TV}}\!}

\def\R{\mathbb{R}}

\def\beq{\begin{equation}}  \def\eeq{\end{equation}}
\def\bb{\begin{eqnarray*}}  \def\ee{\end{eqnarray*}}
\def\b{\begin{eqnarray}}    \def\e{\end{eqnarray}}

\newcommand{\eec}{{\mathrm e}} 	
\newcommand{\RR}{\mathbb{R}}

\newcommand{\norm}[1]{\|#1\|}               			 % norm

\newcommand{\tvnorm}[1]{\norm{#1}_{_{TV}}}                 
 	
\newcommand{\bnorm}[1]{\bigl\|#1\bigr\|_{_{TV}}} % tvnorm with bigger sides
\newcommand{\les}{\leqslant}
\newcommand{\ges}{\geqslant}
\newcommand{\ab}[1]{\vert#1\vert}           			% absolut value
\newcommand{\ii}{{\mathrm i}}               				% complex i
\newcommand{\dd}{{\mathrm{d}}}              			% roman d

\newcommand{\w}{\widehat}
\newcommand{\F}{{\mathcal{F}}} 						% set of all distributions
\newcommand{\ZZ}{\mathbb{Z}}	
\newcommand{\M}{{\cal M}}  								% z- measures
\begin{document}

\title{Multidimensional compound Poisson approximations for symmetric distributions}

\author{ V. \v Cekanavi\v cius\footnote{Corresponding author.}  \, and  S. Jokubauskien\.e \\
 {\small
Institute of Applied Mathematics, Faculty  of Mathematics and Informatics, Vilnius University,}\\
{\small Naugarduko 24, Vilnius 03225, Lithuania.}
\\{\small Department of Mathematics and Statistics, Faculty of Informatics, Vytautas Magnus University,}\\
{\small  Universiteto  10-202, 53361, Akademija, Lithuania}
\\{\small E-mail:
vydas.cekanavicius@mif.vu.lt } {\small and  simona.jokubauskiene@vdu.lt}
 }
\date{}
%\jobname, \today }

\maketitle
%%%%%%%%%%%%%%%%%%%%%%%%%%%%%%%%%%%%%%%%%%%%%%%%%%%%%%%%%%%%%%%%%%%%%%

\begin{abstract}

Distribution of the sum of independent identically distributed symmetric lattice vectors is approximated by the accompanying compound Poisson law and the second-order Hipp-type signed compound Poisson measure. Bergstr\"om-type asymptotic expansion is constructed. The accuracy of approximation is estimated in the total variation metric and, in many cases, is of the order $O(n^{-1})$.

\vspace*{.5cm} \noindent {\emph{Keywords:} \small compound Poisson approximation, convolution of distributions, sums of random vectors, total variation metric, the first uniform Kolmogorov theorem}

\vspace*{.5cm} \noindent {\small {\bf MSC 2000 Subject
Classification}:
Primary 60F99;   % Central limit and other weak theorems
Secondary  60G50; % Sums of independent random variables; random walks
}
\end{abstract}

\newpage

%%%%%%%%%%%%%%%%%%%%%%%%%%%%%%%%%%%%%%%%%%%%%%%%%%%%%%%%%%%%%%%%%%%%%%%%%%%%%%%%%%%%%%%%%%%%%%%%%%%%%%%%%%
\section{Introduction}

 Compound Poisson (CP) approximations have many applications and are one of the popular fields of research in modern probability theory, see, for example,  \cite{BCV01, BUV25, CN24, GMM03, LV22, SU20, SV09} and the references therein.  Research papers in this field are numerous. However, the absolute majority of results are proved for one-dimensional random variables. 
  The situation is very different for random vectors (r.v.s), where results are a few and usually obtained only for Kolmogorov, Prokhorov, or local metrics, see \cite{GZ22, JC25, Z92} and the references therein. Note that, for sequences of r.v.s, a significant part of multidimensional research is devoted to normal approximation \cite{BLX18, BG25} or behavior of tails of distributions \cite{LCW26}. We are primarily interested in  triangular array of vectors.

  %In uniform Kolmogorov metric CP approximations are investigated in \cite{GZ22, Z92}. 
  In this paper,  for sums of lattice identically distributed symmetric r.v.s, we adapt  Le~Cam's convolution approach and properties of metrics to prove results for %a stronger 
  total variation. % metric . 
   Apart from the general case,  when no moment assumptions are imposed, we consider  mixtures of  distributions concentrated on separate lines with finite analogues of variances, demonstrating how  known one-dimensional results can be applied to obtain the multidimensional estimates.

We introduce necessary notation. The sets of all integers and natural numbers are denoted respectively by $\ZZ$ and $\N$.
By $\M_d$  and  $\F_d\subset\M_d$ we denote the set of all finite (signed) $d$-dimensional measures and the set of all $d$-dimensional distributions, respectively. Further on $d<\infty$ is some fixed natural number.
 Let $\,I_{\vec{a}}\in\F_d\,$ be the distribution concentrated at $\vec{a}\in\RR^d$, $I=I_{\vec{0}}$, where $\vec{0}=(0,0,\dots,0)\in\R^d$.  Let $M,V\in\M_d$, then their convolution, for any $d$-dimensional Borel set $A$ is defined as
 $(M*V)\{A\} = \int_{\R^d} M\{ A-\vec{x}\}V\{\dd \vec{x}\}. $

   By $\L(\vec{X})$ we denote the distribution of a r.v. $\vec{X}$. Note that if $S_n=\vec{X}_1+\vec{X}_2+\dots+\vec{X}_n$, where all r.v.s. are independent identically distributed (i.i.d.)  and $\L(\vec{X}_1)=F$, then $\L(S_n)=F^{*n}$, where power is understood in the convolution sense. In this paper we prefer  measure notation over r.v.s, since the main tools used for the proofs are properties of norms and convolutions. 
  
 Exponential measure for $M\in\M_d$ is defined as
$ \exp(M) = \sum_{j\ges 0}{M^{*j}}/{j!}$,
where $M^{*0}=I$. 

Compound Poisson distribution with parameter $\l>0$ and compounding distribution $F\in\F_d$ is defined as  $\exp (\l(F-I))$. Note that CP distribution can be expressed also as
a distribution of random sum $ \vec{Y}=\sum_{j=0}^{\pi_\l} \vec{X}_j$, where $\vec{X},\vec{X}_1,\vec{X}_2,\dots$ are i.i.d. r.v.s  independent from Poisson random variable $\pi_\l$ and   $\L(\,\vec{X})=F$.  CP distribution $\exp (F\!-\!I)$ is called \emph{accompanying distribution} to $F\in\F_d$.  If $F\in\F_d$, but $\l<0$, then $\exp(\l(F\!-\!I))$ is called signed compound Poisson measure (SCP). All CP distributions are infinitely divisible.

 For $F,G\in\F_d$ total variation distance and total variation norm are defined as
\[
\d(F,G):=\sup_B\left|F\left\{B\right\}-G\left\{B\right\}\right|=:\frac{1}{2}\tvnorm{F-G},
  \]
 where supremum is taken over all $d$-dimensional Borel sets. Total variation norm can also be defined via Hahn decomposition, see \cite{CN24}, p. 228. For $F,G\in\F_d$ with supp $F$,  supp  $G\subset \ZZ^d$
 \[\tvnorm{F-G}=\sum_{\vec{m}\in\ZZ^d}\left|F\left\{\vec{m}\right\}-G\left\{\vec{m}\right\}\right|.\]

To make proofs  more concise the same symbol $\,C\,$ is used to denotes all  positive constants.  If ambiguities can arise we supply $C$ with indices.
 Symbol $\,C(N)$ denotes  quantity, which depends only on  $N$. Notation $f(n)=O(n^{-k})$ means that  $f(n)n^k$  is bounded by some $C$ for all $n$. 
Multiplication is superior to division.  For any $x\!\in\!\RR$ let $\lfloor x\rfloor $  denote the integer  part of $x$.

\section{Known results}

General CP approximations for independent symmetric vectors in  Kolmogorov and Prokhorov metrics were investigated in \cite{GZ22,Z24} (see also references, therein). In this paper, we consider stronger total variation metric.  Multidimensional estimates in total variation for CP approximations are few.  Dependent vectors were investigated in Theorem 6.8  \cite{N11}. 
 Non-symmetric vectors were considered by Roos \cite{Roos02,R07, R17}, who used the so-called Kerstan's method. In
   \cite{KC14}, Kerstan's method was applied for symmetric vectors concentrated on coordinate axes of $\ZZ^d$. Set $\bar{e}_1=(1,0,\dots,0)$,..., $\bar{e}_m=(0,\dots,0,1,0,\dots,0)$, 
    $\bar{e}_d=(0,\dots,0,1)$. Let $\vec{X}_1,\dots,\vec{X}_n$ be independent r.v.s  taking values in $\ZZ^d$ and let
for
$m=1,2,\dots,d,\, i=1,2,\dots,n$:
\bb
\L(\vec{X}_i)&=& q_i I+\sum_{m=1}^d p_{i,m}F_m,\quad q_i=\p(\vec{X}_i=\vec{0}), \quad q_i+\sum_{m=1}^d p_{i,m}=1,\label{forma}\\
\hbox{supp}\, F_m&\subset& \left\{k\bar{e}_m,\ k=\pm 1, \pm 2, \dots\right\},\quad F_m\left\{k\bar{e}_m\right\}=F_m\left\{-k\bar{e}_m\right\},\\
\lambda_m&=&\sum_{i=1}^n p_{i,m},\quad\sigma^2_m=2\sum_{k=1}^\infty k^2 F_m\left\{k\bar{e}_m\right\}, \quad
0\les q_i,p_{i,m}\les 1,\label{lambdar}\nonumber\\
\alpha&=&\sum_{i=1}^ng(2(1-q_i))\min\left\{2^{-3/2}\sum_{m=1}^d\frac{p_{i,m}^2}{\lambda_m},\Big(\sum_{m=1}^d
p_{i,m}\Big)^2\right\},\\
g(x)&=&2\sum_{s=2}^\infty\frac{x^{s-2}}{s!}(s-1)=\frac{2\eec^x(\eec^{-x}-1+x)}{x^2},\quad
x\ne 0.
\ee
Observe that in above setting coordinates of $\vec{X}_i$ are uncorrelated and, for example, $\p(\vec{X}_i=(1,1,1,\dots,1)=0$. In \cite{KC14} it was proved that
if  $2\alpha\eec<1$, $\sigma^2_m<\infty$, ($m=1,2,\dots,d$), then

\b\label{KrC1}
\lefteqn{\d \left(\prod_{i=1}^{*n}\L(\vec{X}_i),\prod_{i=1}^{*n}\exp(\L(\vec{X}_i)-I)
\right)}\hskip 2cm\nonumber\\
&&\les
\frac{15.98}{(1-2\alpha\eec)^{3/2}}\sum_{l=1}^d(1+\sigma_l)\sum_{m=1}^d\frac{1}{\lambda_m^2}\sum_{i=1}^n
p_{i,m}^2.
\e

 If $d,\sigma_m,p_{i,m}\asymp C$, then the rate of
accuracy in (\ref{KrC1}) is $O(n^{-1})$.\\

In the case of i.i.d  r.v.s. $\vec{X}_i\equiv \vec{X}$, when $p_{i,m}\equiv p_m$ and $q=\p(\vec{X}=\vec{0})\ges 4/5$,  it was proved that
\b\label{p1is5}
\lefteqn{\d\left(\left(qI+\sum_{m=1}^dp_mF_m\right)^{*n},\quad\exp\left\{n\sum_{m=1}^d p_m(F_m-I)\right\}
\right)}\hskip 4cm\nonumber\\
&&\les
\frac{17.34}{n}\left(\sum_{m=1}^d\sqrt{1+\sigma_m}\right)^2,
\e
see \cite{KC14}, Theorem 3.2. 
Kerstan's method has obvious restrictions. For (\ref{KrC1}) and (\ref{p1is5}) to hold, probabilities  $\p(\vec{X}_i=\vec{0})$ must be large.  Consider, for example, the case $d=1$, $q_i=p_i=1/2$. Then  $2\alpha\eec>1$ and (\ref{KrC1}) cannot be applied. Moreover,  $\sigma_m$ can be large even if supp $\vec{X}_i$ is finite. The assumption that all coordinate vectors are uncorrelated is also very restrictive.\\

%%%%%%%%%%%%%%%%%%%%%%%%%%%%%%%%%%%%%%%%%%%%%%%%
Note that  any distribution with $q=F\{\vec{0}\}>0$ can be written as $F=qI+(1-q)V$. From one-dimensional Poisson approximation to the binomial law it follows then that
\beq\label{binp}
\d\big(F^{*n},\exp(n(F\!-\!I))\big)\les \min(1\!-\!q, n(1\!-\!q)^2),
\eeq
see \cite{CN24}, p. 207 and (3.15) for details.
However, (3) does not take into account the symmetry of $F$ and  requires $q$ to be very close to 1, much closer than needed for (\ref{p1is5}).

\section{Results}

%%%%%%%%%%%%%%%%%%%%%%%%%%%%%%%%%%%%%%%%%%%%%%%%%%%%%%%%%%%%%%%%%%%%%%%%%%%%%%%%%%

\subsection{The general case of vectors with possibly correlated coordinates} 

We will demonstrate that for identically distributed symmetric r.v.s, the  assumptions used for (\ref{p1is5})  can be significantly  weakened. 
Let
\beq\label{CKm1}
\vec{X}\sim F\in\F_d^{s}, \qquad \hbox{supp}\, F\subset \ZZ^d, \qquad q:=F\{\vec{0}\}=\p(\vec{X}=\vec{0})\in (0,1).
\eeq
Assumptions (\ref{CKm1}) are more general than used in (\ref{KrC1}), (\ref{p1is5}). Indeed, in (\ref{p1is5}) it was required that $q\ges 4/5$, meanwhile in (\ref{CKm1}) $q$ can be reasonably  small. Moreover,  non-zero probabilities are not restricted to coordinate axes. Further let $p_{\vec{x}}:=\p(\vec{X}=\vec{x})=F\{\vec{x}\}$ and $\sum_{\vec{x}}$ denotes summation over all $\vec{x}\in \hbox{supp } F\setminus \{\vec{0}\}$. %$\vec{x}\in\ZZ^d\setminus\left\{\vec{0}\right\}$.
 Let
 \beq\label{smix}
 %\delta(y):=2\sum_{\vec{x}}\min\left(p_{\vec{x}},\frac{1}{y\eec }\right).\qquad 
 \delta(y)=\sum_{\vec{x}}\min\left(1,  4y p_{\vec{x}}\right).
  \eeq
 
Note that $\delta(y)\les  4y(1-q)$, $\delta(y/k)\les \delta(y)$ if $k\ges 1$ and $\delta(y)\les 2N$ if  supp $ F=\{\vec{0},\pm\,\vec{x}_1,\dots,\pm\,\vec{x}_{N}\}$.

 The following theorem gives a general idea about what constant can be expected for large $n$.

%%%%%%%%%%%%%%%%%%%%%%%%%%%%%%%%%%%%%%%%%%%%%%%%%%%%%%%%%%%%%%%%%%%%%%%%%%%%%%%%%%%%%%%%%%

\begin{thm}\label{CKMa} Let assumption (\ref{CKm1}) hold. Then, for any $n\in\N$
\beq
\d\big(F^{*n},\exp(n(F\!-\!I))\big)\les \frac{\delta^2(n/2)}{4n}\left(1+\frac{C(\delta(nq)+1)}{q^4 \sqrt{n}}\right)+\frac{C(1-q)^{9/2}}{q^{5}n\sqrt{n}} .\label{cac1}\\
\eeq
\end{thm}

%%%%%%%%%%%%%%%%%%%%%%%%%%%%%%%%%%%%%%%%%%%%%%%%%%%%%%%%%%%%%%%%%%%%%%%%%%%%%%%%%%%%%%
The following shorter version of Theorem \ref{CKMa} with the smaller power of $q$ in denominator holds.

\begin{customthm}{\ref{CKMa}*}\label{star} Let assumption (\ref{CKm1}) hold. Then, for any $n\in\N$
\beq
\d\big(F^{*n},\exp(n(F\!-\!I))\big)\les \frac{C(1+\delta^2(n))}{nq^{7/2}}.\label{cac1star}
\eeq
\end{customthm}

%The following example demonstrates that
 Theorem \ref{CKMa} can be used even if  (\ref{KrC1}) and (\ref{p1is5}) are not applicable. 
%%%%%%%%%%%%%%%%%%%%%%%%%%%%%%%%%%%%%%%%%%%%%%%%%%%%%%%%%%%%%%%%%%%%%%%%%%%%%%%%%%%%%%%%%%%%%%%%%%%%

\begin{exam}\label{1over3}
Let $d=2$ and $\vec{X},\vec{X}_1,\dots,\vec{X}_n$ be i.i.d.  random vectors with 
\bb
&&\p\big(\vec{X}=(0,0)\big)=\frac{1}{2}, \quad \p\big(\vec{X}=(\pm k,0)\big)=\frac{1}{2k(k+1)(k+2)},\qquad k=1,2,\dots\\
&& \p\big(\vec{X}=(0,\pm j)\big)=\frac{1}{2j(j+1)(j+2)},\qquad j=1,2,\dots.\ee
Observe that (\ref{p1is5}) is not applicable, since $q=1/2<4/5$ and  $\sigma_1, \sigma_2=\infty$. On the other hand,
\[
\delta(n)\les 2\left(2\sum_{k=1}^{\lfloor n^{1/3}\rfloor+1}1+n\eec\sum_{k=\lfloor n^{1/3}\rfloor+2}^\infty \frac{1}{k(k+1)(k+2)}\right)<10n^{1/3}\]
and from (\ref{cac1star}) follows estimate
$\d(F^{*n},\exp(n(F\!-\!I)))\les Cn^{-1/3}$.
\end{exam}

Next we use SCP approximation of Hipp's type, see \cite{H86,R07} and \cite{SV09}, Ch. 10.7. Let for any real %(positive or negative)
 $a$
\beq\label{Ds}
D:=\exp\left((F-I)-\frac{(F-I)^{*2}}{2}\right),\quad D^{*a}:=\exp\left(a(F-I)-\frac{a(F-I)^{*2}}{2}\right)
\eeq

%%%%%%%%%%%%%%%%%%%%%%%%%%%%%%%%%%%%%%%%%%%%%%%%

\begin{thm}\label{ZCJ1tv}  Let assumption (\ref{CKm1}) hold. Then, for any $n\in\N$,
\beq\label{cac3}\d(F^{*n},D^{*n})\les \frac{15.75\delta^3(nq/6)}{q^{7/2}\,n^2}\left(1+\frac{C}{q^{7/2}\sqrt{n}}\right)+\frac{C(1\!-\!q)^{15/2}}{q^8\,n^2\sqrt{n}}.\eeq
\end{thm}

%%%%%%%%%%%%%%%%%%%%%%%%%%%%%%%%%%%%%%%%%%%%%%%%
SCP approximation $D$ can be viewed as first order asymptotic expansion in the exponent. Next we consider more usual first order asymptotic expansion to accompanying CP distribution.

\begin{thm}\label{ZCJ2tv} Let assumption (\ref{CKm1}) hold. Then, for any $n\in\N$
\b
\lefteqn{\d\left(F^{*n},\exp(n(F\!-\!I))*\left(I-\frac{n(F\!-\!I)^{*2}}{2}\right)\right)}\hskip 3cm\nonumber\\
&&\les \frac{C(\delta^3(nq)+\delta^4(nq))}{q^{11/2}\, n^2}+\frac{C(1\!-\!q)^6}{q^{13/2}\,n^2}.\label{cac2}
\e
\end{thm}
%%%%%%%%%%%%%%%%%%%%%%%%%%%%%%%%%%%%%%%%%%%%%%%%
%%%%%%%%%%%%%%%%%%%%%%%%%%%%%%%%%%%%%%%%%%%%%%%%
%%%%%%%%%%%%%%%%%%%%%%%%%%%%%%%%%%%%%%%%%%%%%%%%
It is easy to check that for the %case considered in 
Example \ref{1over3} the accuracy in (\ref{cac3})  (respectively (\ref{cac2})) is of the order $O(n^{-1})$ (respectively $O(n^{-2/3})$).

It is known that in Kolmogorov metric the closeness of subsequent convolutions $F^{*n}$ and $F^{*(n+1)}$ can be estimated by $C(d)n^{-1}$, see \cite{Z88}. Assumption (\ref{CKm1}) allows for similar estimate in total variation.

\begin{thm}\label{ZCalphatv}  Let assumption (\ref{CKm1}) hold. Then, for any $n\in\N$
\beq
\d\big(F^{*n},F^{*(n+1)}\big)\les \frac{C(\delta^2(n)+1)}{q^{7/2}\, n}.\label{cac4}\\
\eeq
\end{thm}

Orders of accuracy in (\ref{cac1})--(\ref{cac4}) are explicit, when $\vec{X}$ takes only finite number of values.

\begin{corollary} \label{ZCJcor} Let assumption (\ref{CKm1}) hold and let supp $F=\{\vec{0},\pm\,\vec{x}_1,\pm\,\vec{x}_2,\dots,\pm\,\vec{x}_{N}\}$, $N\ges d$. Then 
\bb
&&\d\big(F^{*n},\exp(n(F\!-\!I))\big)\les \frac{ N^2}{n}\left(1+ \frac{CN}{q^5\sqrt{n}}\right),
\\
&&\d(F^{*n},F^{*(n+1)})\les CN^2q^{-7/2}\,n^{-1},\qquad \d(F^{*n},D^{*n})\les C N^3q^{-7/2}\,n^{-2},\\%\left(1+\frac{C}{q^{9/2}\sqrt{n}}\right),\\
&&\d\left(F^{*n},\exp\big(n(F\!-\!I))*(I-n(F\!-\!I)^{*2}/2)\big)\right)\les CN^4 q^{-13/2}\,n^{-2}.
\ee
\end{corollary}
 Note that only the number of points matter, not their values. The following example shows that even for r.v.s with bounded supports estimate (\ref{cac1}) can be more accurate than  (\ref{p1is5}).
 
 %%%%%%%%%%%%%%%%%%%%%%%%%%%%%%%%%%%%%%%%%%%%%%%%%%%%%%%%%%%%%%%%%%%%%%%%%%%%
 \begin{exam}\label{1isn4} Let $d=2$, $\p(\vec{X}=(\pm 1,0))=\p(\vec{X}=(\pm n,0))=\p(\vec{X}=(0,\pm 1))=\p(\vec{X}=(0,\pm n))=0.025$, $\p(\vec{X}=(0,0))=0.8$. Then  (\ref{p1is5}) gives  trivial order of approximation $O(1)$. Meanwhile  
 from Corollary \ref{ZCJcor} it follows that approximation by $\exp(n(F-I))$ is of the order $O(n^{-1})$.
 \end{exam}
 
 %%%%%%%%%%%%%%%%%%%%%%%%%%%%%%%%%%%%%%%%%%%%%%%%%%%%%%%%%%%%%%%%%%%%%%%%%%
%%%%%%%%%%%
Observe also that estimates in Corollary \ref{ZCJcor} remain meaningful even if $q\to 0$, $N\to\infty$ slowly. For example, if $q=1/\ln n$, $N=O(\ln n)$ the estimates still have quite good order of approximation.

Though absolute constant in Corollary \ref{ZCJcor} is very reasonable, the following result shows that, at least in the scheme of sequences, it might be much smaller. 

\begin{proposition}\label{propo}
Let supp $F=\left\{\vec{0},\pm\,\vec{x}_1,\pm\,\vec{x}_2,\dots,\pm\,\vec{x}_{N}\right\}$, $N\ges d$ and let $N, q, p_{\vec{x}}$ and support of $F$ do not depend on $n$. Then
\beq\label{ascon}
\lim_{n\to\infty}n\, \d(F^{*n},\exp(n(F\!-\!I)))\les 0.17504N+0.05855N(N-1).
\eeq
\end{proposition}

%Let $N,q$ be some fixed constant. Then from  Corollary \ref{ZCJcor} we get
%\[\lim_{n\to\infty}\frac{n}{N^2}\, \d(F^{*n},\exp(n(F\!-\!I)))\les 1,\]
%which gives the idea about possible smallness of positive constants  for large $n$.
 Explicit constants can be obtained for the difference of two CP distributions with the same compounding distribution $F$.

\begin{thm} \label{ZCJexp} Let assumption (\ref{CKm1}) hold and let supp  $F=\{\vec{0},\pm\,\vec{x}_1,\pm\,\vec{x}_2,\dots,\pm\,\vec{x}_N\}$, $N\ges d$. Then, for any $a,b> 0$, the following estimate holds
\[%\label{eabs}
\d\big(\exp(a(F\!-\!I)),\exp(b(F\!-\!I))\big)\les \left(1+\frac{2N+1}{\eec}\right)\frac{\ab{b\!-\!a}}{\max(a,b)}.
\] 
\end{thm}

If $q$ is very large ($F$ is almost degenerate) the following  trivial estimate can be used
\[\d\big(\exp(a(F\!-\!I)),\exp(b(F\!-\!I))\big)\les 2\ab{b\!-\!a}(1\!-\!q).\]

So far we considered short asymptotic expansions. Next we demonstrate that, at least for $D$, it is possible to construct long so-called Bergstr\"om-type \cite{BRG51} asymptotic expansion.
Let 
\beq\label{MK2}
 M_k:=F^{*n}-\sum_{m=0}^k\binom{n}{m}D^{*(n-m)}*(F-D)^{*m}, \qquad k\les n-1.
 \eeq
 
 \begin{thm} \label{BergD} Let assumption (\ref{CKm1}) hold. Then, for any $n\in\N$, $k<n$, the following estimate holds
 \[%\label{bergD2}
 \tvnorm{M_k}\les \frac{C(k)\delta^{3(k+1)}(nq)}{q^{6k+11/2} n^{2(k+1)}} +\frac{C(k)(1-q)^{6(k+1)}}{n^{2(k+1)}q^{6k+13/2}}.
 \]
 \end{thm}

\begin{corollary} \label{Bcor} Let assumption (\ref{CKm1}) hold and let supp $F=\{\vec{0},\pm\,\vec{x}_1,\pm\,\vec{x}_2,\dots,\pm\,\vec{x}_{N}\}$, $N\ges d$. Then 
\[%\label{bergD2}
 \tvnorm{M_k}\les \frac{C(k)N^{3(k+1)}}{q^{6k+13/2} n^{2(k+1)}}.
 \]
\end{corollary}
 
%%%%%%%%%%%%%%%%%%%%%%%%%%%%%%%%%%%%%%%%%%%%%%%%%%%%%%%%%%%%%%%%%%%%%%%%%%%%%%%%%%%%%%%%%%%%%%%%%%%%%%%%%%%%%%
%%%%%%%%%%%%%%%%%%%%%%%%%%%%%%%%%%%%%%%%%%%%%%%%%%%%%%%%%%%%%%%%%%%%%%%%%%%%%%%%%%%%%%%%%%%%%%%%%%%%%%%%%%%%%%%
Theorem \ref{BergD} can be used as intermediate result for construction of longer  expansions. 
%%%%%%%%%%%%%%%%%%%%%%%%%%%%%%%%%%%%%%%%%%%%%%%%%%%%%%%%%%%%%%%%%%%%%%%%%%%%%%%%%%%%%%%%%%%%%%%%%%%%%%%%%%%%%%

\subsection{ The case of r.v.s concentrated on finite number of lines} \label{symspec}

Observe that for some cases estimate (\ref{p1is5}) is more accurate than the one in (\ref{cac1}).

\begin{exam}\label{1isk4} Let $d=2$ and $\vec{X},\vec{X}_1,\dots,\vec{X}_n$ be i.i.d random vectors with  
\bb
&&\p(\vec{X}=(0,0))=\frac{7}{9}, \quad \p(\vec{X}=(\pm k,0))=\frac{1}{2k(k+1)(k+2)(k+3)},,\qquad k=1,2,\dots\\
&& \p(\vec{X}=(0,\pm j))=\frac{1}{2j(j+1)(j+2)(j+3)},,\qquad j=1,2,\dots.\ee
Then, arguing similarly to Example \ref{1over3}, we get
\[
\delta(n)\les 2\left(2\sum_{k=1}^{\lfloor n^{1/4}\rfloor+1}1+2n\eec\sum_{k=\lfloor n^{1/4}\rfloor+2}^\infty \frac{1}{k(k+1)(k+2)(k+3)}\right)\les Cn^{1/4}.\]
From (\ref{cac1star}) it follows that $\d(F^{*n},\exp(n(F\!-\!I)))\les Cn^{-1/2}$. Meanwhile,  $\sigma^2_1,\sigma^2_2\les 2$
and from (\ref{p1is5}) it follows that $\d(F^{*n},\exp(n(F\!-\!I)))\les Cn^{-1}$.
\end{exam}

Below we obtain some generalizations of (\ref{p1is5}). We  assume that random vectors are concentrated on a finite number of lines and satisfy certain moment conditions. 
As before we assume that  $\vec{X},\vec{X}_1,\dots,\vec{X}_n$ are i.i.d. r.v.s in $\RR^d$ , $\L(\vec{X})=F$.
We also assume that
\begin{enumerate}
\item For some natural integer $K\ges d$
\beq\label{sc1}
F=qI+\sum_{m=1}^Kp_mF_m,\quad q,p_1,\dots,p_m\in (0,1),\quad \sum_{m=1}^Kp_m=1-q.
\eeq
\item There exist non-zero vectors  $\vec{y}_1,\vec{y}_2,\dots,\vec{y}_K\in\RR^d$ such that
\beq\label{sc2}
\hbox{supp}\, F_m\subset\{k\vec{y}_m:k\in\ZZ\setminus\{0\}\},\quad \hbox{supp}\, F_m\cap \hbox{supp} \,F_k=\emptyset,\quad k\ne m.
\eeq
\item  All  supp $F_m\ne \emptyset$ and  $F_m$ are symmetric:
\beq\label{sc3}
F_m\big\{k\vec{y}_m\big\}=F_m\big\{-k\vec{y}_m\big\}, \quad k\in\big\{1,2,\dots\big\},\quad m=1,\dots,K.
\eeq
\item %The maximal span of lattice $k\vec{y}_m$ is 1, that is, 
There exist $k^0_m$ such that
\beq\label{sc4}
F_m\big\{k^0_m\vec{y}_m\big\}>0,\quad F_m\big\{(k^0_m+1)\vec{y}_m\big\}>0.
\eeq
\end{enumerate}
Distribution considered in (\ref{p1is5})  is special case of above setting, since it suffices to take  $K=d$ and $\vec{y}_m=\bar{e}_m$. 
We can  define 'projection' of $F$ to $\vec{y}_m$ by considering one-dimensional integer-valued random variable $Y$ such that 
\beq\label{project}
\p(Y=k):=F_m\big\{k\vec{y}_m\big\},\quad k\in\N.\eeq
Then  (\ref{sc4}) is an assumption that $Y$ is concentrated on a lattice with the maximal span equal to one. Though not formulated explicitly, assumption (\ref{sc4}) was used for the proof of (\ref{p1is5}) in \cite{KC14}.
 We  denote  the  variance  of $Y$ by
\[\sigma^2_m=\sum_{k=-\infty}^\infty k^2F_m\{k\vec{y}_m\}=2\sum_{k=1}^\infty k^2F_m\{k\vec{y}_m\}.\]

First we show that assumption $q\ges 4/5$ used in (\ref{p1is5}) can be dropped.

%%%%%%%%%%%%%%%%%%%%%%%%%%%%%%%%%%%%%%%%%%%%%%%%%%%%%%%%%%%%%%%%%%%%%%%%%%%%%%%%%%%%%%%%%%%%%%%%%%%%%%%%%%%%%%%%%%%%%%

\begin{thm}\label{SK1} Let assumptions (\ref{sc1})--(\ref{sc4}) hold. Then, for any $n\in\N$,
\b\label{K1a}
\lefteqn{\d \big(F^{*n}, \exp(n(F\!-\!I))\big)}\hskip 1cm\nonumber\\
&&\les \frac{1.76}{n}\left(\sum_{m=1}^K\sqrt{1+\sigma_m}\right)^2\left(1+\frac{C}{\sqrt{n}\,q^5}\sum_{m=1}^K\sqrt{1+\sigma_m}\right).
\e
\end{thm}

A shorter version with non-explicit constant can also be proved.

\begin{customthm}{\ref{SK1}*}\label{starA} 
Let assumptions (\ref{sc1})--(\ref{sc4}) hold. Then, for any $n\in\N$,
\beq\label{K2}
\d \big(F^{*n}, \exp(n(F\!-\!I))\big)\les \frac{C}{nq^{7/2}}\left(\sum_{m=1}^K\sqrt{1+\sigma_m}\right)^2.
\eeq
\end{customthm}

%%%%%%%%%%%%%%%%%%%%%%%%%%%%%%%%%%%%%%%%%%%%%%%%%%%%%%%%%%%%%%%%%%%%%%%%%%%%%%%%%%%%%%%%%%%%%%%%%%%%%%%%%%%%%%%%%%%%%%%%%%%%%

Next we consider SCP approximation by $D$.

\begin{thm}\label{SKD2a} Let assumptions (\ref{sc1})--(\ref{sc4}) hold. Then, for any $n\in\N$,
\beq\label{K1}
\d (F^{*n}, D^{*n})\les \frac{C}{n^2\,q^{13/2}}\left(\sum_{m=1}^K\sqrt{1+\sigma_m}\right)^3.
\eeq
\end{thm}
%%%%%%%%%%%%%%%%%%%%%%%%%%%%%%%%%%%%%%%%%%%%%%%%%%%%%%%%%%%%%%%%%%%%%%%%%%%%%%%%%%%%%%%%%%%%%%%%%%%%%%%%%%%%%%%%%%%%%%%%%%%%%

We can formulate analogues of Theorems \ref{ZCJ2tv} and \ref{ZCalphatv}.

\begin{thm}\label{SKAs} Let assumptions (\ref{sc1})--(\ref{sc4}) hold. Then, for any $n\in\N$,
\beq
\d\big(F^{*n},\exp(n(F\!-\!I))*(I-n(F\!-\!I)^{*2}/2)\big)\les \frac{C}{n^2\, q^{13/2}}\left(\sum_{m=1}^K\sqrt{1+\sigma_m}\right)^4.\label{cac2d}
\eeq
\end{thm}
%%%%%%%%%%%%%%%%%%%%%%%%%%%%%%%%%%%%%%%%%%%%%%%%%%%%%%%%%%%%%%%%%%%%%%%%%%%%%%%%%%%%%%%%%%%%%%%%%%%%%%%%%%%%%%%%%%%%%%%%%%%%%

\begin{thm}\label{SKAF}  Let assumptions (\ref{sc1})--(\ref{sc4}) hold. Then, for any $n\in\N$,
\[ %\label{K2F}
\d (F^{*n}, F^{*(n+1)})\les \frac{C}{nq^{7/2}}\left(\sum_{m=1}^K\sqrt{1+\sigma_m}\right)^2.
\]
\end{thm}

%%%%%%%%%%%%%%%%%%%%%%%%%%%%%%%%%%%%%%%%%%%%%%%%%%%%%%%%%%%%%%%%%%%%%%%%%%%%%%%%%%%%%%%%%%%%%%%%%%%%%%%%%%%%%%%%%%%%%%%%%%%%%

Finally we formulate analogue of Theorem \ref{ZCJexp}.

\begin{thm}\label{SKAE}  Let assumptions (\ref{sc1})--(\ref{sc4}) hold. Then, for any $a,b>0$,
\[%\label{K2F}
\d \big(\exp(a(F\!-\!I)), \exp(b(F\!-\!I))\big)\les \frac{1.7\ab{b-a}}{\max(a,b)}\sum_{m=1}^K\sqrt{1+\sigma_m}.
\]
\end{thm}

%%%%%%%%%%%%%%%%%%%%%%%%%%%%%%%%%%%%%%%%%%%%%%%%%%%%%%%%%%%%%%%%%%%%%%%%%%%%%%%%%%%%%%%%%%%%%%
Note that by assumption $K\ges d$. Therefore, strictly speaking, estimates in Theorems \ref{SK1} -- \ref{SKAE} depend on $d$ and become trivial if $d\to\infty$.

\subsection{Some numerical simulations}

From the results of this paper it follows that the accuracy of approximation is very good for large $n$. We give some examples of  estimates for simulated data when $n$ is small.
We begin from the case $d=1$. Let $\p(X=-1)=\p(X=1)=p$, $\p(X=0)=1-2p$, $F=\L(X)$, that is
$\w F(t)=q+p(\eec^{-\ii t}+\eec^{\ii t})$.

\begin{table}[t]
\centering
\begin{tabular}{|l|l|l|l|}
  \hline
  % after \\: \hline or \cline{col1-col2} \cline{col3-col4} ...
  $p$ &n=10& n=100 & n=500  \\
  \hline
  0.01 &0.0028088  & 0.0018381&0.0003544 \\
  0.10&0.0182749&0.0017669&0.0003509\\
  0.25&0.0185827& 0.0017586 & 0.0003503 \\
  0.45 &0.0544728&0.0017514 &0.0003503 \\
   \hline
\end{tabular}
\caption{$\d(F^{*n}, \exp(n(F-I)))$  for small $n$.}\label{Tab1}
\end{table} 

Estimates for different simulations ar given in table \ref{Tab1}. On one hand, we see that the accuracy of approximation is very good. On the other hand, if we consider $n\ges 100$, then the accuracy of approximation is $\approx 0.175/n$, which is better than $\approx 1/n$ from Corollary 1. Note also that, for $n=10$ and small $p$ , the accuracy of approximation is determined by (\ref{binp}).

Next let us consider direct extension of the previous simulation to 2-dimensional case.
%\begin{exam}
 Let $q_i=1-2p_i$, $\w F_i(t_i)=q_i+p_i(\eec^{\ii t_i}+\eec^{-\ii t_i})$, $i=1,2$, $\w F(t_1,t_2)=\w F_1(t_1)\w F_2(t_2)$.
%\bb\w F(t_1,t_2)&=& \w F_1(t_1)\w F_2(t_2)
%=q_1q_2+p_1q_2\eec^{\ii t_1}+p_1q_2\eec^{-\ii t_1}
%+q_1p_2\eec^{\ii t_2}+q_1p_2\eec^{-\ii t_2}\\
%&+&p_1p_2\eec^{\ii (t_1+t_2)}+p_1p_2\eec^{\ii (t_1-t_2)}+p_1p_2\eec^{\ii (-t_1+t_2)}
%+p_1p_2\eec^{-\ii(t_1+t_2)}.
%\ee
%Observe that 
%\[\hbox{supp}\ F^{*n}=\left\{-n,-n+1,\dots,n-1,n\right\}\times \left\{-n,-n+1,\dots,n-1,n\right\}.\]
For approximation we use convolution of accompanying distributions
$
\w H(t_1,t_2)=\exp(\w F_1(t_1)-1)\exp(\w F_2(t_2)-1).
$

 In table \ref{Tab2} we present some
simulated estimates  for $\d(F^{*n}, H^{*n})$. From the triangle inequality $\d(F_1^{*n}*F_2^{*n}; \exp(n(F_1-I))*\exp(n(F_2-I)))\les \d(F_1^{*n};\exp(n(F_1-I)))+\d(F_2^{*n};\exp(n(F_2-I)))$ and the first simulation, one can expect the accuracy at least of the order $\approx 0.34/n$, $n\ges 100$. From  table \ref{Tab2} we see  that, in reality, the accuracy is even better  $\approx 0.25/n$.  

\begin{table}[ht]
\centering
\begin{tabular}{|l|l|l|l|l|}
  \hline
  % after \\: \hline or \cline{col1-col2} \cline{col3-col4} ...
  $p_1$ &$p_2$& n=10 & n=100 &n=500 \\
  \hline
  0.01 & 0.01 & 0.0048421& 0.0025247&0.0004952 \\
  0.10&0.10&0.0258073&0.0024872&0.0004974\\
  %0.25& 0.25 & 0.05069 & 0.00498 & 0.00099&0.00049\\
  %0.25 & 0.45 & 0.11152 &0.00499 &0.00099&0.00049 \\
  0.10& 0.45& 0.0552711 & 0.0024914& 0.0004976 \\
  0.45 & 0.45&0.0678198 & 0.0024953 & 0.0004978\\
   \hline
\end{tabular}
\caption{$\d(F^{*n}, H^{*n})$  for small $n$.}\label{Tab2}
\end{table}

%\end{exam}

%%%%%%%%%%%%%%%%%%%%%%%%%%%%%%%%%%%%%%%%%%%%%%%%%%%%%%%%%%%%%%%%%%%%%%%%%%%%%%%%%%%%

\section{Auxiliary results}

First we formulate general facts about exponential measures. Let $M,V\in\M_d$, $k\in\N$. We have %Convolution of two exponential measures is also exponential measure:
$\exp(M)*\exp(V)=\exp(M\!+\!V)$. From definition of exponential measure it follows that    
\beq\label{normB}
\eec^M=I+M+\frac{M^{*2}}{2!}+\cdots+\frac{M^{*k}}{k!}+\frac{M^{*(k+1)}}{k!}*\int_0^1\eec^{\tau M}(1-\tau)^k\dd\tau.
\eeq
 %For the proofs we use the following  properties of total variation norm.
Let $F\in\F_d$, $a>0$. The following well-known relations hold: $\tvnorm{F}=1$,
  \b
&&\tvnorm{M\!*\!V}\les\tvnorm{M}\,\tvnorm{V},\quad \tvnorm{\exp(M)}\les\exp(\tvnorm{M}), \quad  \label{norms}\\ 
%&&,\quad\tvnorm{\exp\{a(F-I)\}}=1, \label{neq1}\\
&&\tvnorm{(F-I)^{*k}*\exp(a(F-I))}\les\tvnorm{(F-I)*\exp(a(F-I)/k)}^k,\label{aka}\\
%&&\eec^M=I+M+\frac{M^{*2}}{2!}+\cdots+\frac{M^{*k}}{k!}+\frac{M^{*(k+1)}}{k!}\int_0^1\eec^{\tau M}(1-\tau)^k\dd\tau, \label{normB}\\
&&\Big\|\sum_{j=0}^\infty \alpha_j F^{*j}\Big\|_{_{TV}}\les %\sum_{j=0}^\infty\ab{\alpha_j}\left\|F^{*j}\right\|_{_{TV}}= 
\sum_{j=0}^\infty\ab{\alpha_j}=
\Big\|\sum_{j=0}^\infty\alpha_j I_j\Big\|_{_{TV}}, \label{cpto1dim}
\e
where  $I_j\in\F_1$ is one-dimensional distribution concentrated at $j$.  % Properties (\ref{norms}) and (\ref{aka}) are well known. %Total variation of any distribution equals 1, therefore (\ref{neq1}) follows.
Sometimes inequality (\ref{cpto1dim}) %  follows from (\ref{norms}), see Eq. (A.3) in \cite{CN24} or Eq. (12) in \cite{Pre85}. it 
is formulated in terms of random sums, see Eq. (3.1), \cite{VeCh96}.
 It allows to reformulate many one-dimensional results for multidimensional case. Next lemma exemplifies  such reformulation. 

%%%%%%%%%%%%%%%%%%%%%%%%%%%%%%%%%%%%%%%%%%%%%%%%%%%%%%%%
\begin{lemma}\label{CN24C.1} If $\l>0$,  $k\in\N$, $\vec{x}\in\RR^d$. Then
 \b 
  && \tvnorm{(I_{-\vec{x}}+I_{\vec{x}}-2I)*\exp\left(\lambda(I_{-\vec{x}}+I_{\vec{x}}-2I)\right)}\les\frac{1}{\l}, \label{oD2}\\
&&\sup_{P\in\F_d}\tvnorm{(P-I)^{*k}\!*\eec^{\l(P-I)}}\les\bigg(\frac{2k}{\eec\l}\bigg)^{k/2}.\label{koD3}
\e
\end{lemma}
   Estimate (\ref{oD2}) follows from Lemma 4.6 in \cite{CeRoAISM06} and (\ref{cpto1dim}).
  Estimate (\ref{koD3}) follows from Proposition A.2.7 in \cite{BHJ}, (\ref{cpto1dim}) and (\ref{aka}).

%%%%%%%%%%%%%%%%%%%%%%%%%%
%%%%%%%%%%%%%%%%%%%%%%%%%%%%%%%%%%%%%%%%%%%%%%%%%%%%%%%%%%%%%%%%%%%%%%%%%%%%%%%%%
The following technical estimates follow  from Eq. (31) in \cite{K86} and  from \cite{roo01} Lemma 6.

\begin{lemma} \label{Kroo} Let $\l>0$. Then
\[
\int_{-\pi}^\pi \sin^2(t/2)\eec^{-2\l \sin^2(t/2)}\dd t\les \frac{\sqrt{\pi}}{\lambda^{3/2}},\qquad
\int_{-\pi}^\pi \eec^{-2\l \sin^2(t/2)}\dd t\les \frac{\pi}{\sqrt{3\l}}\left(1+\sqrt{\frac{\pi}{2}}\right). 
\]
\end{lemma}

%%%%%%%%%%%%%%%%%%%%%%%%%%%%%%%%%%%%%%%%%%%%%%%%%%%%%%%%%%%%%%%%%%%%%%%%%%%%%%%
%%%%%%%%%%%LB.2 generalization
Next we present generalization of Lemma B.2 from \cite{CN24}.

\begin{lemma}\label{B.2new} Let $p_i\in (0,1)$,  ($i=i,2,\dots,n$), $p_0=\max_i p_i$, $\tau\in [0,1]$. Then
\[
\sup_{F\in\F_d}\left\|\exp\left(\frac{(1+p_0)}{2}\sum_{i=1}^np_i(F-I)-\frac{\tau}{2}\sum_{i=1}^np_i^2(F-I)^{*2}\right)\right\|_{_{TV}}\les\frac{3.5}{\sqrt{1-p_0}}.
\]
Particularly, if $p_i\equiv p$, then
\beq\label{e2p}
\sup_{F\in\F_d}\left\|\exp\left(\frac{np(1+p)}{2}(F-I)-\frac{np^2\tau}{2}(F-I)^{*2}\right)\right\|_{_{TV}}\les\frac{3.5}{\sqrt{1-p}}.
\eeq
\end{lemma}

\pr In view of (\ref{cpto1dim}), it  suffices to prove lemma for one-dimensional case $F=I_1\in\F_1$.
Let 
\[M:=\frac{(1+p_0)}{2}\sum_{i=1}^np_i(I_1-I)-\frac{\tau}{2}\sum_{i=1}^np_i^2(I_1-I)^{*2}.\]
We denote by  $\ii$  complex unit, that is, $\ii^2=-1$. Set  $\eec^{\ii x}:=\cos x+\ii\sin x$ and let $\w M(t)$ be Fourier transform of $M$.
Set $a=\sum_{i=1}^n p_i$, $\nu=a(1+p_0)/2$ and
\bb
&&\w M_c(t):=\w M(t)-\ii t\nu=\frac{a(1+p_0)}{2}\left(\eec^{\ii t}-1-\ii t\right)-\frac{\tau}{2}\sum_{i=1}^n p_i^2(\eec^{\ii t}-1)^2,\\
&& b=2.634\sqrt{\frac{a}{1-p_0}}.
\ee

Then
\bb
\left|\eec^{\w M_c(t)}\right|&\les&\exp\left(-(1+p_0)a\sin^2(t/2)+\frac{p_0}{2}a\ab{\eec^{\ii t}-1}^2\right)\\
&\les&\exp\left(-(1-p_0)a\sin^2(t/2)\right).
\ee
Observe that
\[\ab{\w M_c'(t)}\les\frac{a(1+p_0)}{2}\ab{\eec^{\ii t}-1}+ap_0\ab{\eec^{\ii t}-1}\les 4a\ab{\sin(t/2)}.\]
Therefore,
\[\left|\left(\eec^{\w M_c(t)}\right)'\right|^2\les 16 a^2\sin^2(t/2)\eec^{-2(1-p_0)a\sin^2(t/2)}.\]
and  by Lemma \ref{Kroo}
\[
b\int_{-\pi}^\pi\left|\eec^{\w M_c(t)}\right|^2\dd t+\frac{1}{b}\int_{-\pi}^\pi\left|\left(\eec^{\w M_c(t)}\right)'\right|^2\dd t\les 
\frac{21.54}{1-p_0}.
\]
If $b>2.5$, then lemma's statement follows from formula of inversion (see (D.10) in \cite{CN24})
 \[ %\label{TVAP} 
 \tvnorm{\eec^M}^2 \les \frac{1\!+\!b\pi}{2\pi} \int\limits_{-\pi}^{\pi}
 \left( \left| \eec^{\w M_c(t)}\right|^2 +
 b^{-2} \left|\left(\eec^{\w M_c(t)}\right)'\right|^2 \right) \dd t 
 \]
  and a simple estimate
\[\frac{1+b\pi}{2\pi}\les\frac{b}{2}\left(1+\frac{2}{5\pi}\right)\les 0.5637\,b.\]

 If $b\les 2.5$, then $a\les 0.901(1-p_0)$, $a(1+3p_0)\les 0.901(1-p_0)(1+3p_0)\les 1.2014$ and
\[\tvnorm{\eec^M}\les
 \exp\left(\frac{a(1+p_0)}{2}\tvnorm{I_1-I}+\frac{ap_0}{2}\tvnorm{I_1-I}^{2}\right)\les\eec^{a(1+3p_0)}\les 3.33.
\]
\hspace*{\fill}$\Box$  

%%%%%%%%%%%%%%%%%%%%%%%%%%%%%%%%%%%%%%%%%%%%%%%%%%%%%%%%%%%%%%%%%%%%%%%%%%%%%%
%%%%%%%%%%%%%%%%%%%%%%%%%%%%%%%%%%%%%%%%%%%%%%%%%%%%%%%%%%%%%%%%%%%%%%%%%%%%%

\begin{lemma}\label{CeRo4.6} Let $P\in\F_1$ be a symmetric distribution concentrated on $\,\ZZ\!\setminus\!\{0\}.$ Assume that $\,P\,$ has a finite variance $\,\sigma^2\,$. 
Then for any $\l>0$ and
$j\in\N$
\[ %\label{ZRoos1}
 \bnorm{(P-I)^{*j}*\eec^{\l(P-I)}} \les
\frac{3.6\sqrt{1\!+\!\sigma}\,j^{j+1/4} }{\l^{j}\eec^{j}}
\]
\end{lemma}

Lemma \ref{CeRo4.6} is  part of Lemma 4.6 from \cite{CeRoAISM06}.
%%%%%%%%%%%%%%%%%%%%%%%%%%%%%%%%%%%%%%%%%%%%%%%%%%%%%%%%%%%%%%%%%%%%%%%%%%%%%%%
Similar bounds are valid for the compound binomial distribution. Next lemma is Lemma 4 in \cite{Roo00}.
 \begin{lemma}\label{CN24C} Let $\,q\!=\!1\!-\!p.$ If $\,0<p< 1$, $\,j,n\in \N$, then 
 \b
\sup_{F\in\F_d}  \tvnorm{(F\!-\!I)^{*j}*(qI+pF)^{*n}} %\!&\le&\! 
	%\binom{n\!+\!j}{j}^{-1/2}(pq)^{-j/2}\nonumber \\ 
   &\le&\!\sqrt{\eec}\,j^{1/4} \bigg(\frac{n}{n\!+\!j}\bigg)^{n/2}
	\bigg( \frac{j}{(n\!+\!j)pq}\bigg)^{j/2}.\label{C6a}% \label{a5bin}
    \e
\end{lemma}

We recall that notation $\Theta$ is used for all measures satisfying $\tvnorm{\Theta}\les 1$ .
%%%%%%%%%%%%%%%%%%%%%%%%%%%%%%%%%%%%%%%%%%%%%%%%%%%%%%%%%%%%%%%%%%%%%%%%%%
\begin{lemma}\label{D2irF}  Let $F=qI+pV\in\F_d$, $q=1-p\in(0,1)$, $m\in\N$. Then
for any $m\in\N$, $a\in\RR$
\b
&&D-I=C(F-I)*\Theta, \quad \tvnorm{D^{*a}}\les\eec^{4\ab{a}},\label{FD00}\\
&&D^{*m}=\frac{3.5}{\sqrt{q}}\exp\left(\frac{m q(F-I)}{2}\right)*\Theta_m,\quad \bnorm{D^{*m}}\les\frac{3.5}{\sqrt{q}},\label{D2Ftrys}\\
%&&F-D=\eec^4 (F-I)^{*3}*\Theta=\eec^4 p^3(V-I)^{*3}*\Theta,\label{FD2}\\
&&F-D=\frac{(F-I)^{*3}}{3}+C_1(F-I)^{*4}*\Theta=
C_2(F-I)^{*3}*\Theta\nonumber\\
&&=C_2p^3(V-I)^{*3}*\Theta.
\label{FD2}
\e
\end{lemma}

\pr The first estimate in (\ref{D2Ftrys}) follows from (\ref{e2p}). The second estimate in (\ref{D2Ftrys}) follows from the fact that the total variation norm of any CP distribution is equal to 1 and , therefore,
\[ \bnorm{D^{*m}}\les\frac{3.5}{\sqrt{q}}\bnorm{\exp\left(m q(F-I)/2\right)}\tvnorm{\Theta_m}\les \frac{3.5}{\sqrt{q}}.\]
Observe that $\tvnorm{F-I}\les\tvnorm{F}+\tvnorm{I}= 2$ and, therefore, by (\ref{norms}) 
%\[\tvnorm{D}\les\exp\{\tvnorm{F-I}+\tvnorm{F-I}^2/2\}\les \eec^4.\]
%Similarly
\[\left\|D^{*a}\right\|_{_{TV}}
%&:=&\left\|\exp\left\{a(F-I)-a(F-I)^{*2}/2\right\}\right\|_{_{TV}}\\
%&\les&
\les\exp(\ab{a}\tvnorm{F-I}+\ab{a}\tvnorm{F-I}^2/2)\les\eec^{4\ab{a}}.
\]
Expressions (\ref{FD00}) and (\ref{FD2}) follow directly from definition of exponential measure.
\hspace*{\fill}$\Box$  \\

%%%%%%%%%%%%%%%%%%%%%%%%%%%%%%%%%%%%%%%%%%%%%%%%%%%%%%%%%%%%%%%%%%

The following result is of special interest, since it shows that long asymptotic expansions can be constructed without assumptions of symmetry. Let $M_k$ be defined by (\ref{MK2}).
%\[ M_k:=F^{*n}-D^{*n}-\sum_{m=1}^k\binom{n}{m}D^{*(n-m)}*(F-D)^{*m}.\]

\begin{lemma} \label{m3} Let $F=qI+pV\in\F_d$, $q=1-p\in(0,1)$, $k,n\in\N$, $k\les n-1$. Then
\[\tvnorm{M_k}\les \frac{C(k)p^{3(k+1)/2}}{n^{(k+1)/2} q^{(3k+4)/2}}.\]
\end{lemma}

\pr
    If $n\les 12$, then
\[\tvnorm{M_k}\les 1+\sum_{m=0}^k\binom{n}{m}\eec^{4(n-m)}(1+\eec^{4})^m\les 1+(1+2\eec^4)^{12}= C=\frac{n^kC}{n^k}\les\frac{12^kC}{n^k}.\]
Therefore,  we further assume $n\ges 13$. 
Bergstr\"om identity  \cite{BRG51} allows to write
\[M_k=\sum_{m=k+1}^n\binom{m-1}{k}F^{*(n-m)}*(F-D)^{*(k+1)}*D^{*(m-k-1)}.\]
Noting that
 \beq \label{bergidentA}
 \sum_{m=k+1}^n\binom{m-1}{k}=\binom{n}{k+1}
  \eeq 
and applying  (\ref{FD2}), (\ref{D2Ftrys}), (\ref{koD3}), (\ref{C6a}), (\ref{bergidentA})  we obtain
\bb
\tvnorm{M_k}&\les&\sum_{m=k+1}^n\binom{m-1}{k}\bnorm{F^{*(n-m)}*(F-D)^{*(k+1)}*D^{*(m-1)}}\\
&\les&C(k)\sum_{m\les \lfloor n/2\rfloor}\binom{m-1}{k}\bnorm{(F-I)^{*3(k+1)}*F^{*(n-m)}}\bnorm{D^{*(m-k-1)}}\\
&+&C(k)\sum_{m>\lfloor n/2\rfloor}\binom{m-1}{k}\bnorm{F^{*(n-m)}}\bnorm{(F-I)^{*3(k+1)}*D^{*(m-k-1)}}\\
&\les& \frac{C(k)n^{k+1}}{\sqrt{q}}p^{3(k+1)}\bnorm{(V-I)^{*3(k+1)}*(qI+pV)^{*\lfloor n/2\rfloor}}\\
&+& \frac{C(k)}{\sqrt{q}}n^{k+1} p^{3(k+1)}\left\| (V-I)^{*3(k+1)}*\eec^{npq(V-I)/6}\right\|_{_{TV}}\\
&\les& \frac{C(k)p^{3(k+1)/2}}{n^{(k+1)/2} q^{(3k+4)/2}}. \ee
\hspace*{\fill}$\Box$  
%%%%%%%%%%%%%%%%%%%%%%%%%%%%%%%%%%%%%%%%%%%%%%%%%%%%%%%%%%%%%%%%%%%%%%%%%%%%%%%%%%%%%%%%%%%%%%%%%%%%%%%%%%%%%%

\begin{lemma}\label{Fexp} Let assumptions (\ref{CKm1}) hold. Then, for any $a>0$, $k\in\N$, the following estimate holds
\beq\label{fexp}
\tvnorm{(F-I)^{*k}*\exp(a(F-I))}\les \left(\frac{k}{2 a}\right)^k\delta^k(a/k)\les \frac{C(k)\delta^k(a)}{ a^k},
\eeq
where $\delta(a)$ is defined by (\ref{smix}).\end{lemma}

\pr  Since $F$ is symmetric, we can write
\beq\label{simi3}
F-I=\sum_{\vec{x}}p_{\vec{x}}\left(I_{\vec{x}}-I\right)=\sum_{\vec{x}}p_{\vec{x}}\left(I_{-\vec{x}}-I\right)=\frac{1}{2}\sum_{\vec{x}}p_{\vec{x}}\left(I_{\vec{x}}+I_{-\vec{x}}-2I\right).
\eeq
Total variation of CP distribution equals 1. Therefore,
\[\exp(a(F-I))=\exp\left(p_{\vec{x}}\left(I_{\vec{x}}+I_{-\vec{x}}-2I\right)\right)
*\Theta,
\]
for any $\vec{x}$ from  $F$  support and some $\Theta=\Theta(F)$ such that $\tvnorm{\Theta}\les 1$.

Therefore, applying (\ref{oD2}) we get
\bb%\label{simi4}
\lefteqn{\tvnorm{(F\!-\!I)*\exp(a(F\!-\!I))}}\hskip 2cm\\
&&\les\frac{1}{2}\sum_{\vec{x}}p_{\vec{x}}\left\|(I_{\vec{x}}\!+\!I_{-\vec{x}}\!-\!2I)*\exp\left(p_{\vec{x}}\left(I_{\vec{x}}+I_{-\vec{x}}-2I\right)\right)\right\|_{_{TV}}\\
&&\les \frac{1}{2}\sum_{\vec{x}}p_{\vec{x}}\min\left(4,\frac{1}{a  p_{\vec{x}}}\right)=\frac{1}{2a}\delta(a).
\ee
Estimate for $k>1$ follows from (\ref{aka}).
\hspace*{\fill}$\Box$  

%%%%%%%%%%%%%%%%%%%%%%%%%%%%%%%%%%%%%%%%%%%%%%%%%%%%%%%%%%%%%%%%%%%%%%%%%%%%%%%%%%%%%%%%%%%%%%%%%%%%%

\begin{lemma}\label{Fdexp} Let assumptions (\ref{sc1})--(\ref{sc4}) hold. Then, for any $a>0$, $k\in\N$, the following estimate holds
\beq\label{dexp}
\tvnorm{(F-I)^{*k}*\exp(a(F-I))}\les \frac{(3.6k)^k}{\eec^k a^k}\left(\sum_{m=1}^K\sqrt{1+\sigma_m}\right)^k.
\eeq
\end{lemma}

\pr Observe that $\eec^{a(F-I)}$ can be expressed as convolution of two CP distributions 
\[\eec^{a(F-I)}=\exp(ap_m(F_m-I))*\Theta,\quad\tvnorm{\Theta}=\bigg\|\exp\bigg(a\sum_{j\ne m}^K p_j(F_j-I)\bigg)\bigg\|_{_{TV}}=1,\]
where the last equality follows from the fact that total variation norm of any distribution equals 1.
 Noting that estimates for convolutions of $F_m$ can be treated as estimates for one-dimensional distributions (see (\ref{project})) and applying  Lemma \ref{CeRo4.6} we get
\bb
\tvnorm{(F-I)*\exp(a(F-I))}&\les& \sum_{m=1}^Kp_m \left\|(F_m-I)*\exp(ap_m(F_m-I))\right\|_{_{TV}}
\\
&\les&\frac{3.6}{a\eec}\sum_{m=1}^K\sqrt{1+\sigma_m}.
\ee
For the proof, when $k>1$, we use (\ref{aka}).\hspace*{\fill}$\Box$  

%%%%%%%%%%%%%%%%%%%%%%%%%%%%%%%%%%%%%%%%%%%%%%%%%%%%%%%%%%%%%%%%%%%%%%%%%%%%%%%%%%%%%%%%%%%%%%%%%%%%%

%We recall that $\Theta$ is used for various measures, satisfying $\tvnorm{\Theta}\les 1$.

\begin{lemma}\label{D2accF} Let $F=qI+pV\in\F_d$, $q\in(0,1)$, and let $D$ be defined by (\ref{Ds}). Then for any $n\in\N$ the following equalities hold
\b
D^{*n}-\eec^{n(F-I)}\!&=&\!\frac{3.5n}{2\sqrt{q}}(F\!-\!I)^{*2}*\eec^{nq(F-I)/2}*\Theta,\label{D2F.1}\\
D^{*n}-\eec^{n(F-I)}\!&=&\!-\frac{n}{2}\eec^{n(F-I)}\!*\!(F\!-\!I)^{*2}+\frac{3.5n^2}{4\sqrt{q}}\eec^{nq(F-I)/2}\!*\!(F\!-\!I)^{*4}\!*\!\Theta.\label{D2F.2}
\e
\end{lemma}

\pr Observe that $F-I=p(V-I)$. Then applying (\ref{normB}) and (\ref{e2p}) we obtain
\bb
\lefteqn{D^{*n}-\eec^{n(F-I)}=\eec^{n(F-I)}*\left(\eec^{-n(F-I)^{*2}/2}-I\right)}\hskip 1cm\\
&&=-\frac{n}{2}\,\eec^{n(F-I)}*(F-I)^{*2}*\int_0^1
\eec^{-n\tau(F-I)^{*2}/2}\dd\tau\\
&&=\frac{n}{2}\,\eec^{nq(F-I)/2}*(F-I)^{*2}*\int_0^1(-1)\eec^{np(1+p)(V-I)/2-\tau np^2(V-I)^{*2}/2}\dd\tau\\
&&=\frac{3.5n}{2\sqrt{q}}\,\eec^{nq(F-I)/2}*(F-I)^{*2}*\Theta,
\ee
which proves (\ref{D2F.1}). The proof of (\ref{D2F.2}) is similar:

\bb
\lefteqn{D^{*n}-\eec^{n(F-I)}}\hskip 0.9cm\\
&&=-\frac{n}{2}\,\eec^{n(F-I)}*(F-I)^{*2}+\eec^{n(F-I)}*\left(\eec^{-n(F-I)^{*2}/2}-I+n(F-I)^{*2}/2\right)\\
%&&=\eec^{nq(F-I)/2}*\frac{n^2}{4}(F-I)^{*4}*\int_0^1\eec^{np(1+p)(V-I)/2-\tau np^2(V-I)^{*2}/2}(1\!-\!\tau)\dd\tau\\
%&&-\frac{n}{2}\eec^{n(F-I)}*(F-I)^{*2}\\
&&=-\frac{n}{2}\,\eec^{n(F-I)}*(F-I)^{*2}+\frac{3.5n^2}{4\sqrt{q}}\eec^{nq(F-I)/2}*(F-I)^{*4}*\Theta.
\ee
\hspace*{\fill}$\Box$  
%%%%%%%%%%%%%%%%%%%%%%%%%%%%%%%%%%%%%%%%%%%%%%%%%%%%%%%%%%%%%%%%%%%%%%%%%%%%%%%%%%%%%%%%%%%%%%%%%%%%%

%%%%%%%%%%%%%%%%%%%%%%%%%%%%%%%%%%%%%%%%%%%%%%%%%%%%%%%%%%%%%%%%%%%%%%%%%%%%%%%%%%%%%%%%%%%%%% 
\section{Proofs}

%%%%%%%%%%%%%%%%%%%%%%%%%%%%%%%%%%%%%%%%%%%%%%%%%%%%%%%%%%%%%%%%%%%%%%%%%%%%%%%%%%%%%%%%%%%%%%%%%%%%%%%%%%%%%%%%%%%%%%%%%%%%%%%%%%
%%%%%%%%%%%%%%%%%%%%%%%%%%%%%%%%%%%%%%%%%%%%%%%%%%%%%%%%%%%%%%%%%%%%%%%%%%%%%%%%%%%%%%%%%%%%%%%%%%%%%%%%%%%%%%%%%%%%%%%%%%%%%%%%%%%
\textbf{Proofs of Theorem \ref{CKMa} and  \ref{SK1}.}  We recall that total variation distance is half of the total variation norm. 
From Lemma \ref{m3} it follows that
\b\label{ACC1}
\left\|F^{*n}-\eec^{n(F-I)}\right\|_{_{TV}}&\les&\frac{C(1-q)^{9/2}}{n\sqrt{n}\,q^{5}}+\left\|D^{*n}-\eec^{n(F-I)}\right\|_{_{TV}}\nonumber\\
&+&\sum_{m=1}^2\binom{n}{m}\left\|D^{*(n-m)}*(F-D)^{*m}\right\|_{_{TV}}.
\e

Taking into account  (\ref{D2Ftrys}), (\ref{FD2}) and (\ref{koD3}) and observing that $\tvnorm{F-I}\les 2$  we obtain
\bb%\label{ACC3}
\lefteqn{\sum_{m=1}^2\binom{n}{m}\left\|D^{*(n-m)}*(F-D)^{*m}\right\|_{_{TV}}\les
C\sum_{m=1}^2n^m\tvnorm{D^{-m}}\left\|D^{*n}*(F-I)^{*3m}\right\|_{_{TV}}}\hskip 0.2cm\nonumber\\
&&\les\frac{C}{\sqrt{q}}\sum_{m=1}^2n^m\left\|
\eec^{nq(F-I)/2}*(F-I)^{*3m}\right\|_{_{TV}}   %\nonumber\\
%&&
\les \frac{C\sqrt{n}}{\sqrt{q}}\left\|
\eec^{nq(F-I)/4}*(F-I)^{*2}\right\|_{_{TV}}\nonumber\\
&&\times\sum_{m=1}^2 n^{m-1/2} \left\|
\eec^{nq(F-I)/4}*(F-I)^{*(2m-1)}\right\|_{_{TV}}\tvnorm{F-I}^{m-1}\nonumber\ee
\bb
&&\les \frac{C\sqrt{n}}{q^2}\left\|
\eec^{nq(F-I)/4}*(F-I)^{*2}\right\|_{_{TV}}.
%\frac{C\delta^2(nq)}{q^4 n\sqrt{n}}.
\ee
Applying  (\ref{fexp})  we get
\beq\label{ACC4}
\sum_{m=1}^2\binom{n}{m}\left\|D^{*(n-m)}*(F-D)^{*m}\right\|_{_{TV}}\les \frac{C\delta^2(nq/4)}{n\sqrt{n}\, q^{4}}\les \frac{C\delta^2(n/2)}{n\sqrt{n}\, q^{4}}.
\eeq 
Applying (\ref{dexp}) we get
\beq\label{ACC5}
\sum_{m=1}^2\binom{n}{m}\left\|D^{*(n-m)}*(F-D)^{*m}\right\|_{_{TV}}\les \frac{C}{n\sqrt{n}\, q^{4}}\left(\sum_{m=1}^K\sqrt{1+\sigma_m}\right)^2.
\eeq

From (\ref{D2F.2}) and (\ref{koD3}) it follows that
\bb%\label{ACC6}
\lefteqn{
\tvnorm{D^{*n}-\eec^{n(F-I)}}\les\frac{n}{2}\left\|(F-I)^{*2}*\eec^{n(F-I)}\right\|_{_{TV}}+\frac{3.5n^2}{4\sqrt{q}}\left\|(F-I)^{*4}*\eec^{nq(F-I)/2}\right\|_{_{TV}}}\hskip 1.9cm\nonumber\\
&&\les\frac{n}{2}\left\|(F-I)^{*2}*\eec^{n(F-I)}\right\|_{_{TV}}\\
&&+\frac{3.5n^2}{4\sqrt{q}}\left\|(F-I)^{*3}*\eec^{nq(F-I)/4}\right\|_{_{TV}}
\left\|(F-I)*\eec^{nq(F-I)/4}\right\|_{_{TV}}\nonumber\\
&&\les\frac{n}{2}\left\|(F-I)^{*2}*\eec^{n(F-I)}\right\|_{_{TV}}+\frac{Cn\sqrt{n}}{q}\left\|(F-I)^{*3}*\eec^{nq(F-I)/4}\right\|_{_{TV}}\!.
\ee
Applying  (\ref{fexp})  we get
\beq\label{ACC7}
\tvnorm{D^{*n}-\eec^{n(F-I)}}\les\frac{\delta^2(n/2)}{2n}+\frac{C\delta^3(nq/4)}{n\sqrt{n}\,q^4}
\les\frac{\delta^2(n/2)}{2n}+\frac{C\delta^2(n/2)\delta(nq)}{n\sqrt{n}\,q^4}.
\eeq
Applying  (\ref{dexp})  we get
\beq\label{ACC8}
\tvnorm{D^{*n}-\eec^{n(F-I)}}\les\frac{25.92}{\eec^2\,n}\left(\sum_{m=1}^K\sqrt{1+\sigma_m}\right)^2+\frac{C}{n\sqrt{n}\,q^4}\left(\sum_{m=1}^K\sqrt{1+\sigma_m}\right)^3.
\eeq

Substituting  (\ref{ACC4}) and (\ref{ACC7}) into (\ref{ACC1}) we complete the proof of (\ref{cac1}).
Substituting  (\ref{ACC5}) and (\ref{ACC8}) into (\ref{ACC1}) we complete the proof of (\ref{K1a}). \hspace*{\fill}$\Box$\\
%%%%%%%%%%%%%%%%%%%%%%%%%%%%%%%%%%%%%%%%%%%%%%%%%%%%%%%%%%%%%%%%%%%%%%%%%%%%%%%%%%%%%%%%%%%%%%%%%%%%%%%%%%%%%%%%%%%%%%%%%%%%%%%%%%
%%%%%%%%%%%%%%%%%%%%%%%%%%%%%%%%%%%%%%%%%%%%%%%%%%%%%%%%%%%%%%%%%%%%%%%%%%%%%%%%%%%%%%%%%%%%%%%%%%%%%%%%%%%%%%%%%%%%%%%%%%%%%%%%%%
%%%%%%%%%%%%%%%%%%%%%%%%%%%%%%%%%%%%%%%%%%%%%%%%%%%%%%%%%%%%%%%%%%%%%%%%%%%%%%%%%%%%%%%%%%%%%%%%%%%%%%%%%%%%%%%%%%%%%%%%%%%%%%%%%%%
%%%%%%%%%%%%%%%%%%%%%%%%%%%%%%%%%%%%%%%%%%%%%%%%%%%%%%%%%%%%%%%%%%%%%%%%%%%%%%%%%%%%%%%%%%%%%%%%%%%%%%%%%%%%%%%%%%%%%%%%%%%%%%%%%%%

\noindent\textbf{Proofs of Theorem \ref{star} and  Theorem \ref{starA}.} Applying Lemma \ref{m3}, (\ref{D2F.1}), (\ref{D2Ftrys}), (\ref{FD2}) we get
\bb\label{BCC1}
\lefteqn{\left\|F^{*n}-\eec^{n(F-I)}\right\|_{_{TV}}\les\tvnorm{M_1}+\left\|D^{*n}-\eec^{n(F-I)}\right\|_{_{TV}}+n\left\|D^{*(n-1)}*(F-D)\right\|_{_{TV}}}\hskip 1cm\nonumber\\
&&\les \frac{C(1\!-\!q)^3}{nq^{7/2}}+\frac{C n}{\sqrt{q}}\left\|(F-I)^{*2}*\eec^{nq(F-I)/4}\right\|_{_{TV}}\left\|\eec^{nq(F-I)/4}\right\|_{_{TV}}\nonumber\\
&&+\frac{C n}{\sqrt{q}}\tvnorm{D^{*(-1)}}\left\|(F-I)^{*2}*\eec^{nq(F-I)/4}\right\|_{_{TV}}\tvnorm{F-I}\Big\|\eec^{nq(F-I)/4}\Big\|_{_{TV}}\nonumber\\
&&\les\frac{C(1\!-\!q)^3}{nq^{7/2}}+\frac{Cn}{\sqrt{q}}\left\|(F-I)^{*2}*\eec^{nq(F-I)/4}\right\|_{_{TV}}.
\ee
It remains to apply (\ref{fexp}) or (\ref{dexp}). \hspace*{\fill}$\Box$\\

%%%%%%%%%%%%%%%%%%%%%%%%%%%%%%%%%%%%%%%%%%%%%%%%%%%%%%%%%%%%%%%%%%%%%%%%%%%%%%%%%%%%%%%%%%%%%%%%%%%%%%%%%%%%%%%%%%%%%%%%%%%%%%%%%%%%%%%%%
%%%%%%%%%%%%%%%%%%%%%%%%%%%%%%%%%%%%%%%%%%%%%%%%%%%%%%%%%%%%%%%%%%%%%%%%%%%%%%%%%%%%%%%%%%%%%%%%%%%%%%%%%%%%%%%%%%%%%%%%%%%%%%%%%%%

\noindent\textbf{Proof of Theorem \ref{ZCJ1tv}.}  From Lemma \ref{m3} it follows that
\beq\label{D2pr1}
\tvnorm{F^{*n}-D^{*n}}\les%\tvnorm{M_4}\les 
\frac{C(1-q)^{15/2}}{\,q^{8}\,n^2\sqrt{n}}+\sum_{m=1}^4\binom{n}{m}\left\|D^{*(n-m)}*(F-D)^{*m}\right\|_{_{TV}}.
\eeq

Therefore, applying (\ref{D2Ftrys}), (\ref{FD2}), (\ref{oD2}) and (\ref{fexp}) %and \ref{B.2new}
 we prove that
\b\label{D2pr2}
\lefteqn{\sum_{m=2}^4\binom{n}{m}\left\|D^{*(n-m)}*(F-D)^{*m}\right\|_{_{TV}}\les\sum_{m=2}^4\binom{n}{m}\left\|D^{*n}*(F-D)^{*m}\right\|_{_{TV}}\eec^{4m}}\hskip 0.5cm\nonumber\\
%&&\les \frac{C}{\sqrt{q}}\sum_{m=2}^4\binom{n}{m}\left\|\exp\left\{\frac{nq(F-I)}{2}\right\}*(F-D)^{*m}\right\|_{_{TV}}\\
&&\les\frac{C}{\sqrt{q}}\sum_{m=2}^4n^m\bnorm{\eec^{nq(F-I)/2}*(F-I)^{*3m}}
\les\frac{C}{\sqrt{q}}\bnorm{\eec^{nq(F-I)/4}*(F-I)^{*3}}\nonumber\\
&&\times\sum_{m=2}^4n^m\bnorm{\eec^{nq(F-I)/4}*(F-I)^{*(2m-1)}}\tvnorm{F-I}^{m-2}\nonumber\\
&&\les \frac{C}{\sqrt{q}}\left\|\eec^{nq(F-I)/12}*(F-I)\right\|_{_{TV}}^3\sum_{m=2}^4\frac{n^m 2^{m-2}}{(nq)^{m-1/2}}\les\frac{C\delta^3(nq/12)}{q^7 n^2\sqrt{n}}\nonumber\\
&&\les\frac{C\delta^3(nq/6)}{q^7 n^2\sqrt{n}}.
\e
Similarly applying  (\ref{FD2}), (\ref{FD00}),  (\ref{D2Ftrys}), (\ref{koD3}) and (\ref{fexp}) we get
\b\label{D2pr3}
\lefteqn{n\tvnorm{D^{*(n-1)}*(F-D)}\les\frac{n}{3}\tvnorm{D^{*(n-1)}*(F\!-\!I)^{*3}}+
Cn\tvnorm{D^{*(n-1)}*(F\!-\!I)^{*4}}
}\hskip 0.3cm\nonumber\\
&&\les \frac{n}{3}\tvnorm{D^{*n}*(F\!-\!I)^{*3}}+\frac{n}{3}\tvnorm{D^{*(-1)}}\tvnorm{D^{*n}*(I-D)*(F\!-\!I)^{*3}}\nonumber\\
&&+ Cn\tvnorm{D^{*(-1)}}\tvnorm{D^{*n}*(F\!-\!I)^{*4}}\nonumber\\
&&\les \frac{3.5n}{3\sqrt{q}}\left\|\eec^{nq(F\!-\!I)/2}*(F\!-\!I)^{*3}\right\|_{_{TV}}+\frac{Cn}{\sqrt{q}}
\left\|\eec^{nq(F-I)/2}*(F\!-\!I)^{*4}\right\|_{_{TV}}\nonumber\\
&&\les \frac{31.5\,\delta^3(nq/6)}{q^{7/2}n^2} +\frac{Cn}{\sqrt{q}}\left\|\eec^{nq(F-I)/4}*(F\!-\!I)\right\|_{_{TV}}
\left\|\eec^{nq(F-I)/4}*(F-I)^{*3}\right\|_{_{TV}}
\nonumber\\
&&\les \frac{31.5\,\delta^3(nq/6)}{q^{7/2}n^2} +\frac{C\sqrt{n}}{q}
\left\|\eec^{nq(F-I)/4}*(F\!-\!I)^{*3}\right\|_{_{TV}}\nonumber\\
&&\les\frac{31.5\,\delta^3(nq/6)}{q^{7/2}n^2}+\frac{C\delta^3(nq)}{q^4 \, n^2\sqrt{n}}.\e
Collecting estimates (\ref{D2pr1}), (\ref{D2pr2}) and (\ref{D2pr3}) we complete the proof of Theorem. \hspace*{\fill}$\Box$\\

%%%%%%%%%%%%%%%%%%%%%%%%%%%%%%%%%%%%%%%%%%%%%%%%%%%%%%%%%%%%%%%%%%%%%%%%%%%%%%%%%%%%%%%%%%%%%%%%%%%%%%%%%%%%%%%%%%%%%%%%%%%%%%%%%%%%%

\noindent\textbf{Proofs of Theorems \ref{ZCJ2tv} and  \ref{SKAs}.}
From Lemma \ref{m3} it follows that
\b\label{ACA1}
\lefteqn{\left\|F^{*n}-\eec^{n(F-I)}*\left(I-\frac{n}{2}(F-I)^{*2}\right)\right\|_{_{TV}}}\hskip 2cm\nonumber\\
&&\les \frac{C(1\!-\!q)^6}{n^2q^{13/2}}+\sum_{m=1}^3\binom{n}{m}\left\|D^{*(n-m)}*(F-D)^{*m}\right\|_{_{TV}}\nonumber\\
&&+\left\|D^{*n}-\eec^{n(F-I)}*\left(I-n(F-I)^{*2}/2\right)\right\|_{_{TV}}
\e

Taking into account (\ref{D2Ftrys}), (\ref{FD2}) and (\ref{koD3})  we obtain
\b\label{ACA2}
\lefteqn{\sum_{m=1}^3\binom{n}{m}\left\|D^{*(n-m)}*(F\!-\!D)^{*m}\right\|_{_{TV}}\les
C\sum_{m=1}^3\binom{n}{m}\tvnorm{D^{*(-m)}}\left\|D^{*n}*(F\!-\!I)^{*3m}\right\|_{_{TV}}
}\hskip 3.45cm\nonumber\\
&&\les
\frac{C}{\sqrt{q}}\sum_{m=1}^3\binom{n}{m}\eec^{4m}\left\|\eec^{nq(F-I)/2}*(F\!-\!I)^{*3m}\right\|_{_{TV}}\nonumber\\
&&\les
\frac{Cn}{\sqrt{q}}\left\|
\eec^{nq(F\!-\!I)/4}*(F\!-\!I)^{*3}\right\|_{_{TV}}\nonumber\\
&&\times\sum_{m=1}^3 n^{m-1} \left\|
\eec^{nq(F-I)/4}*(F\!-\!I)^{*2(m-1)}\right\|_{_{TV}}\tvnorm{F\!-\!I}^{m-1}\nonumber\\
&&\les \frac{Cn}{\sqrt{q}}\left\|
\eec^{nq(F-I)/4}*(F\!-\!I)^{*3}\right\|_{_{TV}}\sum_{m=2}^3\frac{n^{m-1}2^{m-1}}{n^{m-1}q^{m-1}}\nonumber\\
&&\les \frac{Cn}{q^{5/2}}
\left\|
\eec^{nq(F-I)/4}
*(F\!-\!I)^{*3}
\right\|_{_{TV}}. 
\e
From (\ref{D2F.2}) it follows that
\beq\label{ACA4}
\left\|D^{*n}-\eec^{n(F-I)}*\left(I-\frac{n}{2}(F-I)^{*2}\right)\right\|_{_{TV}}\les \frac{Cn^2}{\sqrt{q}}\left\|\eec^{nq(F-I)/2}*(F-I)^{*4}\right\|_{_{TV}}.
\eeq

To proof  (\ref{cac2}) (respectively (\ref{cac2d})) it suffices in (\ref{ACA2}) -- (\ref{ACA4}) to use (\ref{fexp}) (respectively (\ref{dexp}))  and substitute resulting estimates  into (\ref{ACA1}).
 \hspace*{\fill}$\Box$\\

%%%%%%%%%%%%%%%%%%%%%%%%%%%%%%%%%%%%%%%%%%%%%%%%%%%%%%%%%%%%%%%%%%%%%%%%%%%%%%%%%%%%%%%%%%%%%%%%%%%%%%

\noindent\textbf{Proofs of Theorems \ref{ZCalphatv} and \ref{SKAF}.} By the triangle inequality
\bb\left\|F^{*n}-F^{*(n+1)}\right\|_{_{TV}}&\les&\left\|F^{*n}-\eec^{n(F-I)}\right\|_{_{TV}}+\left\|\eec^{n(F-I)}-\eec^{(n+1)(F-I)}\right\|_{_{TV}}\\
&+&
\left\|F^{*(n+1)}-\eec^{(n+1)(F-I)}\right\|_{_{TV}}.\ee
From definition of exponential measure it follows that
\[
\left\|\eec^{n(F-I)}-\eec^{(n+1)(F-I)}\right\|_{_{TV}}=\left\|\eec^{n(F-I)}*\left(I-\eec^{F-I}\right)\right\|_{_{TV}}
\les C\left\|\eec^{n(F-I)}*(F-I)\right\|_{_{TV}}.
\]
It remains to apply (\ref{cac1star}) and (\ref{fexp}) or (\ref{K2}) and (\ref{dexp}).
\hspace*{\fill}$\Box$\\

%%%%%%%%%%%%%%%%%%%%%%%%%%%%%%%%%%%%%%%%%%%%%%%%%%%%%%%%%%%%%%%%%%%%%%%%%%%%%%%%%%%%%%%%%%%%%%%%%%%%%%
%%%%%%%%%%%%%%%%%%%%%%%%%%%%%%%%%%%%%%%%%%%%%%%%%%%%%%%%%%%%%%%%%%%%%%%%%%%%%%%%%%%%%%%%%%%%%%%%%%%%%%%%
\noindent\textbf{Proof of Theorem \ref{ZCJexp}.}
Without loss of generality let $a<b$. From definition of exponential measure (\ref{normB}), the fact that total variation of any distribution equals 1 and Lemma \ref{Fexp} we obtain
\b\label{ADA1}
\lefteqn{\left\|\eec^{b(F-I)}-\eec^{a(F-I)}\right\|_{_{TV}}=\left\|\eec^{a(F-I)}*\left(\eec^{(b-a)(F-I)}-I\right)\right\|_{_{TV}}}\hskip 2.6cm\nonumber\\
&&=\left\|(b-a)\eec^{a(F-I)}*(F-I)*\int_0^1\eec^{\tau(b-a)(F-I)}\dd\tau\right\|_{_{TV}}\nonumber\\
&&\les(b-a)\left\|\eec^{a(F-I)}*(F-I)\right\|_{_{TV}}\int_0^1\left\|\eec^{\tau(b-a)(F-I)}\right\|_{_{TV}}\dd\tau\nonumber\\
&&\les(b-a)\left\|\eec^{a(F-I)}*(F-I)\right\|_{_{TV}}\les \frac{2(b-a)\delta(a)}{a\eec}\nonumber\\
&&\les\frac{2(b-a)(2N+1)}{a\eec}.
\e
If $b\les a(1+\eec/(2N+1))$, then 
\[\frac{2(b-a)(2N+1)}{a\eec}\les \frac{2(b-a)}{b}\left(1+\frac{2N+1}{\eec}\right)\]
and theorem's statement follows from (\ref{ADA1}).
If $b> a(1+\eec/(2N+1))$, then directly
\[\left\|\eec^{b(F-I)}-\eec^{a(F-I)}\right\|_{_{TV}}\les\bnorm{\eec^{b(F-I)}}+\bnorm{\eec^{a(F-I)}}= 2\les \frac{2(b\!-\!a)}{b}\left(1+\frac{2N\!+\!1}{\eec}\right).\]
Recall that total variation distance is half of the total variation norm. \hspace*{\fill}$\Box$\\

%%%%%%%%%%%%%%%%%%%%%%%%%%%%%%%%%%%%%%%%%%%%%%%%%%%%%%%%%%%%%%%%%%%%%%%%%%%%%%%%%%%%%%%%%%%%%%%%%%%%%%%%%%%%%%%%%%%%%%%%%%%%%%%%%%%%%%%%%%%%%
\noindent\textbf{Proof of Theorem \ref{BergD}.}
Applying  Lemma \ref{m3}, (\ref{D2Ftrys}), (\ref{koD3}), (\ref{FD2}) and (\ref{fexp})  we  obtain
\bb
\tvnorm{M_s}&\les&\tvnorm{M_{4s+3}}+\sum_{m=s+1}^{4s+3}\binom{n}{m}\left\|D^{*(n-m)}*(F-D)^{*m}\right\|_{_{TV}}\\
&\les& \frac{C(s)(1\!-\!q)^{6(s+1)}}{n^{2(s+1)}q^{6s+13/2}}+C\sum_{m=s+1}^{4s+3}n^m\tvnorm{D^{*(-m)}}\left\|D^{*n}*(F-I)^{*3m}\right\|_{_{TV}}\\
&\les& \frac{C(s)(1\!-\!q)^{6(s+1)}}{n^{2(s+1)}q^{6s+13/2}}+\frac{C(s)}{\sqrt{q}}\sum_{m=s+1}^{4s+3}n^m\left\|\eec^{nq(F-I)/2}*(F-I)^{*3m}\right\|_{_{TV}}\\
&\les&\frac{C(s)(1\!-\!q)^{6(s+1)}}{n^{2(s+1)}q^{6s+13/2}}+\frac{C(s) n^{s+1}}{\sqrt{q}}\left\|\eec^{nq(F-I)/4}*(F-I)^{*3(s+1)}\right\|_{_{TV}}\\
&\times&\left(1+\sum_{m=1}^{3s+2}n^m\left\|\eec^{nq(F-I)/4}*(F-I)^{*2m}\right\|_{_{TV}}\tvnorm{F-I}^m\right)\\
&\les&\frac{C(s)(1\!-\!q)^{6(s+1)}}{n^{2(s+1)}q^{6s+13/2}}+\frac{C(s)\delta^{3(s+1)}(nq)}{q^{6s+11/2}\,n^{2(s+1)}}.
\ee
\hspace*{\fill}$\Box$\\

%%%%%%%%%%%%%%%%%%%%%%%%%%%%%%%%%%%%%%%%%%%%%%%%%%%%%%%%%%%%%%%%%%%%%%%%%%%%%%%%%%%%%% 
\noindent\textbf{Proof of Theorem \ref{SKD2a}.}
Applying Lemma \ref{m3}, (\ref{D2Ftrys}), (\ref{FD2}), (\ref{koD3}) and (\ref{dexp}) we get
\bb
\tvnorm{F^{*n}-D^{*n}}&\les&\tvnorm{M_3}+C\sum_{m=1}^3\binom{n}{m}\tvnorm{D^{*(-m)}}\tvnorm{D^{*n}*(F-I)^{*3m}}\\
&\les&\frac{C(1\!-\!q)^6}{n^2q^{13/2}}+\frac{C}{\sqrt{q}}\sum_{m=1}^3\eec^{4m}n^m\left\|\eec^{nq(F-I)/2}*(F-I)^{*3m}\right\|_{_{TV}}\\
&\les&\frac{C(1\!-\!q)^6}{n^2q^{13/2}}+\frac{Cn}{\sqrt{q}}\left\|\eec^{nq(F-I)/4}*(F-I)^{*3}\right\|_{_{TV}}\\
&\times&
\sum_{m=1}^3n^{m-1}\left\|\eec^{nq(F-I)/4}*(F-I)^{*2(m-1)}\right\|_{_{TV}}2^{m-1}\\
&\les&\frac{C(1\!-\!q)^6}{n^2q^{13/2}}+\frac{C}{n^2 q^{7/2}}\left(\sum_{m=1}^K\sqrt{1+\sigma_m}\right)^3\sum_{m=2}^3\frac{n^{m-1} 2^{m-1}}{n^{m-1} q^{m-1}}\\
&\les& \frac{C}{n^2\,q^{13/2}}\left(\sum_{m=1}^K\sqrt{1+\sigma_m}\right)^3.
\ee
 \hspace*{\fill}$\Box$\\

%%%%%%%%%%%%%%%%%%%%%%%%%%%%%%%%%%%%%%%%%%%%%%%%%%%%%%%%%%%%%%%%%%%%%%%%%%%%%%%%%%%%%%%%%%%%%%%%%%%%%%%%
\noindent\textbf{Proof of Theorem \ref{SKAE}.}
Without loss of generality let $a<b$. Repeating the proof of (\ref{ADA1}) but using (\ref{dexp}) we prove that
\[%\label{SADA1}
\left\|\eec^{b(F-I)}-\eec^{a(F-I)}\right\|_{_{TV}}\les(b-a)\left\|\eec^{a(F-I)}*(F-I)\right\|_{_{TV}}\les \frac{3.6(b-a)}{a\eec}\sum_{m=1}^K\sqrt{1+\sigma_m}.\]
If $b\les 2.5a$, then 
\[\frac{3.6(b-a)}{a\eec}\sum_{m=1}^K\sqrt{1+\sigma_m}\les \frac{3.4(b-a)}{b}\sum_{m=1}^K\sqrt{1+\sigma_m}. \]
and theorem's statement follows from analogue of (\ref{ADA1}).
If $b> 2.5a$, then directly
\[\bnorm{\eec^{b(F-I)}-\eec^{a(F-I)}}\les 2\les \frac{3.4(b\!-\!a)}{b}\les \frac{3.4(b\!-\!a)}{b}\sum_{m=1}^K\sqrt{1\!+\!\sigma_m}.\]
 \hspace*{\fill}$\Box$

%%%%%%%%%%%%%%%%%%%%%%%%%%%%%%%%%%%%%%%%%%%%%%%%%%%%%%%%%%%%%%%%%%%%%%%%%%%%%%%%%%%%%%%%%%%%%%%%%%%%%%%%%%%%%%
%%%%%%%%%%%%%%%%%%%%%%%%%%%%%%%%%%%%%%%%%%%%%%%%%%%%%%%%%%%%%%%%%%%%%%%%%%%%%%%%%%%%%%%%%%%%%%%%%%%%%%%%%%%%%%
%%%%%%%%%%%%%%%%%%%%%%%%%%%%%%%%%%%%%%%%%%%%%%%%%%%%%%%%%%%%%%%%%%%%%%%%%%%%%%%%%%%%%%%%%%%%%%%%%%%%%%%%%%%%%%
%%%%%%%%%%%%%%%%%%%%%%%%%%%%%%%%%%%%%%%%%%%%%%%%%%%%%%%%%%%%%%%%%%%%%%%%%%%%%%%%%%%%%%%%%%%%%%%%%%%%%%%%%%%%%%
\noindent\textbf{Remarks on simulations.} R-studio, library(cubature) and formulas of inversion were used. Example of the program is given below
 \begin{lstlisting}[language=R]
>n<-10; p1<-0.10; p2<-0.10; q1<-1-2*p1; q2<-1-2*p2; S<-0;  Z<-0; U<-0; for (m in -n:n){integrandz1<- function(x) (q1+2*p1*cos(x[1]))^n*cos(m*x[1])/(2*pi); z1= hcubature(integrandz1, lowerLimit=-pi, upperLimit=pi);  integranda1<- function(x) exp(2*n*p1*(cos(x[1])-1))*cos(x[1]*m)/(2*pi); a1= hcubature(integranda1,lowerLimit=-pi,upperLimit=pi);
 for (k in -n:n){integrandz2<- function(x) (q2+2*p2*cos(x[1]))^n*cos(k*x[1])/(2*pi); z2= hcubature(integrandz2, lowerLimit=-pi, upperLimit=pi); integranda2<- function(x) exp(2*n*p2*(cos(x[1])-1))*cos(x[1]*k)/(2*pi); a2= hcubature(integranda2, lowerLimit=-pi, upperLimit=pi); S=S+abs(z1$integral*z2$integral-a1$integral*a2$integral); Z=Z+ a1$integral*a2$integral; U<-U+z1$integral*z2$integral }}; print(U); print(S); print(Z); var<-S+1-Z;  dtv<-var/2; print(dtv)
 \end{lstlisting}

%% Table %%
%%%%%%%%%%%%%%%%%%%%%%
%\begin{table}
%\caption{}\label{}
%\end{table}

%% Figure %%
%%%%%%%%%%%%%%%%%%%%%%
%\begin{figure}[t]
%\includegraphics{}
%\caption{}\label{}
%\end{figure}

%% Theorem %%
%%%%%%%%%%%%%%%%%%%%%%
%%\begin{thm}[]\label{} Theorem
%%\end{thm}

%% Proof %%
%%%%%%%%%%%%%%%%%%%%%%
%\begin{proof}
%\end{proof}

%% Appendices %%
%%%%%%%%%%%%%%%%%%%%%%
%\begin{appendix}
%\section{Appendix section}
%\end{appendix}

 %%Acknowledgements %%
%%%%%%%%%%%%%%%%%%%%%%
\section*{Acknowledgement}
We are grateful to the referees for their valuable remarks that allowed to improve the initial version of the paper.

%% Funding %%
%%%%%%%%%%%%%%%%%%%%%%
%\begin{funding}
%\end{funding}

%% Bibliography %%
%%%%%%%%%%%%%%%%%%%%%%
%\bibliographystyle{bib/vmsta2-mathphys}
%\bibliography{bib/biblio}

%\bibliography{bib/paperCJ}

{\small \end{document}